\documentclass[a4paper,11pt]{article}
\usepackage{latexsym}
\usepackage{amsmath}
\usepackage{amssymb}
\usepackage{enumerate}
\usepackage{theorem}
\usepackage{array}
\newtheorem{itheorem}{Theorem}

\newtheorem{theorem}{Theorem}[section]
\newtheorem{lemma}[theorem]{Lemma}
\newtheorem{corollary}[theorem]{Corollary}
\newtheorem{proposition}[theorem]{Proposition}
\theorembodyfont{\rmfamily}
\newtheorem{definition}[theorem]{Definition}
\newtheorem{notation}[theorem]{Notation}

\newtheorem{remark}[theorem]{Remark}

\newtheorem{explanation}[theorem]{}
\newcommand{\proof}{\noindent \mbox{\em Proof.\hspace*{2mm}}}
\newcommand{\qed}{\hfill \mbox{$  \Box $}}
\newcommand{\gyokan}{\vskip 4pt}
\title{Rational curves with many rational points over a finite field}
\author{Satoru Fukasawa
\thanks{Partially supported by Grant-in-Aid
for Young Scientists (22740001), JSPS.}
\\
 Department of Mathematical Sciences,
Yamagata University\\
Kojirakawa-machi 1-4-12\\
Yamagata 990-8560, Japan\\
s.fukasawa@sci.kj.yamagata-u.ac.jp
\and 
Masaaki Homma
\thanks{Partially supported by Grant-in-Aid
for Scientific Research (21540051), JSPS.}
\\
 Department of Mathematics,
Kanagawa University\\
Yokohama 221-8686, Japan\\
homma@kanagawa-u.ac.jp
\and
Seon Jeong Kim
\thanks{Partially supported by Basic Science Research Program through the National Research Foundation of Korea(NRF) funded by the Ministry of Education, Science and Technology (2010-0028027).}\\
 Department of Mathematics and RINS\\
Gyeongsang National University\\
Jinju 660-701, Korea \\
skim@gnu.kr
}
\date{}

\begin{document}
\maketitle
\begin{abstract}
We study a particular plane curve over a finite field whose normalization is
of genus $0$.
The number of rational points of this curve achieves the Aubry-Perret bound for rational curves.
The configuration of its rational points and a generalization of the curve
are also presented.
\\
{\em MSC} (2010):14G15, 14G05, 14G50, 14H50
\\
{\em Key Words}: Rational curve, Finite field, Rational point
\end{abstract}

\section{Introduction}
Let $C$ be a curve of degree $d$ in projective plane ${\Bbb P}^2$
over a finite field ${\Bbb F}_q$.
We are interested in the number $N_q(C)$ of the set $C({\Bbb F}_q)$ of 
${\Bbb F}_q$-points.
Especially we want to give a good upper bound for $N_q(C)$ in terms of
$d$ and $q$ for curves $C$ in a certain class.
In 2010, the second and the third authors proved a fact of
this direction \cite{hom-kim2010}.
\begin{itheorem}
If a plane curve $C$ of degree $d \geq 2$ over ${\Bbb F}_q$
has no ${\Bbb F}_q$-linear components, then
\begin{equation}\label{SziklaiBound}
N_q(C) \leq (d-1)q + 1
\end{equation}
unless $d=q=4$ and $C$ is projectively equivalent to
\[
K: (X+Y+Z)^4+(XY+YZ+ZX)^2+XYZ(X+Y+Z)=0
\]
over ${\Bbb F}_4$.
In this exceptional case, $N_4(K) =14$.
\end{itheorem}
Since the upper bound (\ref{SziklaiBound}) was originally
conjectured by Sziklai \cite{szi},
we refer it as the Sziklai bound.
Note that the bound (\ref{SziklaiBound}) makes sense only if
$2 \leq d \leq q+2$ because it is worse than
the obvious bound:
\[
N_q(C) \leq {}^{\#} {\Bbb P}^2 ({\Bbb F}_q) = q^2 + q +1
\]
if $d> q+2$.

The study of curves that attain the Sziklai bound is still under way
\cite{hom-kim2009, hom-kim2011},
but we have not yet met any example of a curve with singularities
which attains it.
We guess that there will be a better bound for curves with singularities,
which is our motivation for focusing on them.

Let $C'$ be an irreducible curve of degree $d$ in ${\Bbb P}^2$
over ${\Bbb F}_q$ whose normalization is ${\Bbb P}^1$.
Since the morphism ${\Bbb P}^1 \to C'$ given by the normalization
is defined over ${\Bbb F}_q$, the number of nonsingular ${\Bbb F}_q$-points
of $C'$ is at most $q+1$. 
Therefore
\begin{equation}\label{AubryPerretBound}
N_q(C') \leq q+1 + \frac{1}{2}(d-1)(d-2),
\end{equation}
because the number of singularities of $C'$ is
at most $\frac{1}{2}(d-1)(d-2)$.
This bound is a special case of
a result of Aubry and Perret \cite[Prop. 2.3]{aub-per}.
So we refer this bound as the Aubry-Perret bound for rational curves.
The bound (\ref{AubryPerretBound})  is, of course,
better than (\ref{SziklaiBound}) in the meaningful range of $d$.

Let $B$ be the rational plane curve over ${\Bbb F}_q$ defined by
the image of
\begin{equation}\label{BallicoHefezCurve}
{\Bbb P}^1 \ni (s,t) \to (s^{q+1}, s^{q}t+st^{q}, t^{q+1}) \in {\Bbb P}^2.
\end{equation}
$B$ is of degree $q+1$ and $N_q(B)$ actually attains the Aubry-Perret bound
for $d=q+1$.
This curve over an algebraically closed field containing ${\Bbb F}_q$
appeared in Ballico and Hefez's classification list \cite[Th. 1]{bal-hef}
of non-reflexive plane curves of degree $q+1$ with second-order $q$
in a different parametrization
\[
(s,t) \to (t^{q+1}, s(s-t)^q, s^q(s-t)) \in {\Bbb P}^2
\]
from (\ref{BallicoHefezCurve}).
Actually those parametrizations are equivalent over ${\Bbb F}_{q^2}$,
but not over ${\Bbb F}_q$.
However we refer the curve $B$ parametrized by (\ref{BallicoHefezCurve})
as the Ballico-Hefez curve.
Recently, in \cite{fuk} the first author has studied the Galois points of $B$
by using a parametrization
\[
(s,t) \to (s^{q+1}, (s+t)^{q+1}, t^{q+1}) \in {\Bbb P}^2
\]
equivalent to (\ref{BallicoHefezCurve}) over
${\Bbb F}_q$, and has found that the constellation of Galois points of $B$
is described in a similar way that the second author did in \cite{hom} for Hermitian curves.

As Hermitian curves have many interesting properties in finite geometry
including coding theory,
Ballico-Hefez curves also might have lovely properties
because of this similarity\footnote{
For other similarity, see Remark~\ref{otherSim}.
}.

In Section~2, we verify that the number of ${\Bbb F}_q$-points of $B$
is actually $q+1+\frac{1}{2}q(q-1)$,
in other words $B$ is a singular maximal curve, and show that
the zeta function of $B$ is given by
$
\frac{(1+T)^{\frac{q^2-q}{2}}}{(1-T)(1-qT)}.
$
In Section~3, we consider the case $q$ is odd.
We give combinatorial characterizations of
the set of ${\Bbb F}_q$-points of $B$,
and compute the $\tau_i$'s,
where $\tau_i$ is the number of ${\Bbb F}_q$-lines
that are $i$-secants to $B({\Bbb F}_q)$.
In Section~4, we compute the $\tau_i$'s for $q$ even.
By using those results, we compute parameters of codes coming from
$B({\Bbb F}_q)$ in Section~5.
Some of them have the largest minimum distance
under fixed length and dimension.
In the last section,
we propose a generalization of the Ballico-Hefez curve,
which is a rational curve in ${\Bbb P}^n$
parametrized by elementary symmetric polynomials
in $t, t^q, \dots , t^{q^{n-1}}$.
We give a formula of the number of ${\Bbb F}_q$-points of
this curve.

\section{Arithmetic of the curve $B$}
We study the arithmetic properties of the Ballico-Hefez curve $B$.
We prepare some notations.
$\Phi$ denotes the morphism (\ref{BallicoHefezCurve}), that is,
$\Phi(s,t) = (s^{q+1}, s^{q}t+st^{q}, t^{q+1})$, and $\varphi$ denotes
$\Phi | {\Bbb A}^1$, where ${\Bbb A}^1= \{ s \neq 0 \}$, that is,
$\varphi (t) = (1, t + t^q , t^{q+1})$.
For a point $P=(\alpha , \beta ) \in {\Bbb P}^1$,
$P^q$ denotes the point $(\alpha^{q} , \beta^{q} ) \in {\Bbb P}^1$.
${\rm Sing}\, B$ denotes the set of singularities of $B$, and
${\rm Flex}\, B$ the set of inflection points of $B$.
The precise definition of an ``inflection point" in our case
will be given in \ref{Order} below.
We also take interest in the order-sequence and
the $q$-Frobenius order-sequence for $\Phi$.
The notion of $q$-Frobenius order-sequence was introduced by
St\"{o}hr and Voloch \cite{sto-vol} as a tool to bound the number of
${\Bbb F}_q$-points of curves.
Here we give the definition of them in a little more general setting than
our original one.

\begin{explanation}\label{Order}
{\em Order-sequence and $q$-Frobenius order-sequence.}
Let $f:\tilde{C} \to {\Bbb P}^n$ be a morphism from a nonsingular curve $\tilde{C}$
over $\overline{\Bbb F}_q$,
where $\overline{\Bbb F}_q$ denotes the algebraic closure of ${\Bbb F}_q$.
In a neighborhood of an assigned point $P \in \tilde{C}$,
$f$ can be represented by $(n+1)$-tuple regular functions
$(f_0, \ldots, f_n)$ one of which is the constant function $1$.
Let $t$ be a local parameter at $P$. Then regular functions around $P$
can be embedded into the formal power series ring
$\overline{\Bbb F}_q[[t]]$ via the identification with the completion of
the local ring at $P \in \tilde{C}$.
The $i$-th Hasse derivation $D^{(i)}$ on $\overline{\Bbb F}_q[[t]]$
is given by $D^{(i)}t^k = \binom{k}{i}t^{k-i}$.
Then $\{ (\varepsilon_0, \dots , \varepsilon_n) 
     \mid \varepsilon_0 < \dots < \varepsilon_n, \ 
     \det \left( (D^{(\varepsilon_i)}f_j)(P)\right)_{i,j} \neq 0 \}$
is nonempty.
The minimum $(n+1)$-tuple in the above set by the lexicographic order
is called the Hermitian $P$-invariant for $f$.
The Hermitian $P$-invariant is constant if $P$ is in a certain nonempty
open subset of $\tilde{C}$.
The Hermitian $P$-invariant of a point in this open subset is called
the order-sequence for $f$.
For details, consult \cite[Ch.7]{hir-kor-tor}.

Let $C = f(\tilde{C}) \subset {\Bbb P}^n$, and
$\varepsilon_0 < \dots < \varepsilon_n$ the order-sequence for $f$.
For a nonsingular point $P' = f(P) \in C$,
the linear subspace spanned by
$\nu +1$ vectors
$\{((D^{(\varepsilon_i)}f_0)(P), \ldots , (D^{(\varepsilon_i)}f_n)(P) )
 | i = 0, 1, \ldots , \nu\}$
in ${\Bbb P}^n$ is the tangent $\nu$-plane at $P'$,
which is denoted by $T_{P'}^{(\nu)}C$.

If the image of $P'$ by $q$-Frobenius map lies on $T_{P'}^{(\nu)}C$
with some $\nu <n$ for almost all $P' \in C$,
the curve $C$ is said to be $q$-Frobenius nonclassical, and
the minimum number $\nu$ having the above property is called
the $q$-Frobenius index of $C$ (see, \cite[Prop. 2]{gar-hom}).
Let $\nu$ be the $q$-Frobenius index of $C$.
The sequence
$\{ \varepsilon_0, \varepsilon_1, \dots , \varepsilon_n \} 
    \setminus \{\varepsilon_{\nu} \}$
is called the $q$-Frobenius order-sequence of $C$.

For a plane curve $C =f(\tilde{C}) \subset {\Bbb P}^2$,
a nonsingular point $P' \in C$ is an inflection point
if $i(C.T_{P'}^{(1)}(C); P')> \varepsilon_2$.
\end{explanation}

Now we go back to our original setting.
\begin{theorem}\label{theorem1}
\begin{enumerate}[{\rm (i)}]
\item For any $P \in {\Bbb P}^1$, the map induced by $\Phi$ on tangent spaces
\[
d{\Phi}_P: T_{P, {\Bbb P}^1} \to T_{\Phi(P), {\Bbb P}^2}
\]
is injective.
\item ${\rm Sing}\, B$ consists of $\frac{q^2 - q}{2}$ ordinary double points.
\item For $P \in {\Bbb P}^1$, $\Phi(P) \in {\rm Sing}\, B$
if and only if
$P \in {\Bbb P}^1({\Bbb F}_{q^2}) \setminus {\Bbb P}^1({\Bbb F}_{q})$.
In this case, $\Phi^{-1}(\Phi(P)) = \{ P, P^q \}$.
\item The order-sequence for $\Phi$ is $\{0, 1, q\}$,
and $B$ is $q$-Frobenius nonclassical.
\item $\Phi(P)$ is an inflection point of $B$ if and only if
$P \in {\Bbb P}^1({\Bbb F}_{q})$.
\item $B({\Bbb F}_{q}) = {\rm Flex}\, B \cup {\rm Sing}\, B$.
\end{enumerate}
\end{theorem}
\proof
The assertion (i) is obvious because
\begin{equation}\label{dPhi}
 \left(
  \begin{array}{c}
   \frac{\partial\Phi}{\partial s}\\
   \frac{\partial\Phi}{\partial t}
  \end{array}
 \right)
    =
 \left(
   \begin{array}{ccc}
     s^q & t^q & 0 \\
     0   & s^q & t^q
   \end{array}
 \right).
\end{equation}
By (i), $\Phi(P)$ is a singular point if and only if
${}^{\#}\Phi^{-1}(\Phi(P)) >1$.
We want to determine such points.
For $P_{\infty} = (0, 1)$, obviously
$\Phi^{-1}(\Phi(P_{\infty})) = \{ P_{\infty} \}$.
For two distinct points $P_{\alpha} = (1, \alpha)$
and $P_{\beta}= (1, \beta) \in {\Bbb A}^1$,
$\Phi(P_{\alpha}) = \Phi(P_{\beta})$ if and only if
$\varphi(\alpha) = \varphi(\beta)$, that is,
\[
\left\{
  \begin{array}{ccc}
   \alpha + \alpha^q &=& \beta + \beta^q \\
   \alpha^{q+1}      &=& \beta^{q+1}.
  \end{array}
\right. 
\]
Hence $\Phi(P_{\alpha}) = \Phi(P_{\beta})$
implies that
$
(X-\alpha)(X-\alpha^q)=(X-\beta)(X-\beta^q).
$
Since $\alpha \neq \beta$,
we have $\alpha = \beta^q$ and $\alpha^q=\beta$.
Hence $P_{\alpha}, P_{\beta} 
\in {\Bbb P}^1({\Bbb F}_{q^2}) \setminus {\Bbb P}^1({\Bbb F}_{q})$
with $P_{\alpha}^q = P_{\beta}$.
Conversely, this condition obviously
leads to $\Phi(P_{\alpha}) = \Phi(P_{\beta})$.
Therefore
$\Phi^{-1}({\rm Sing}\, B) = 
{\Bbb P}^1({\Bbb F}_{q^2}) \setminus {\Bbb P}^1({\Bbb F}_{q})$
and
${\Bbb P}^1({\Bbb F}_{q^2}) \setminus {\Bbb P}^1({\Bbb F}_{q})
\stackrel{\Phi}{\to} {\rm Sing}\, B$
is a $2$-to-$1$ map.
Moreover, since the tangent line to $B$ at the branch corresponding
$P_{\alpha} \in {\Bbb P}^1$ is
\begin{equation}\label{tangent}
 \left|
   \begin{array}{ccc}
    X & Y       & Z \\
    1 &\alpha^q & 0 \\
    0 &   1     &\alpha^q
   \end{array}
 \right|= 0,
\end{equation}
$\Phi(P_{\alpha}) = \Phi(P_{\alpha^q}) \in {\rm Sing}\, B$
is an ordinary double point, that is,
those two branches have different tangent lines.
So (ii) and (iii) have been established, and furthermore
we have also established
\begin{equation}\label{ratptB}
B({\Bbb F}_{q})= \Phi({\Bbb P}^1({\Bbb F}_{q})) \cup {\rm Sing}\, B
\end{equation}
because $\alpha + \alpha^q$ and $\alpha^{q+1}$ are the trace and the norm
from ${\Bbb F}_{q^2}$ to ${\Bbb F}_{q}$ respectively.

Next we calculate the order-sequence for $\Phi$.
Since $d\Phi \neq 0$, the first two orders are $0$ and $1$.
Let $D_t^{(\nu)}$ be the $\nu$-th Hasse derivation on the function field
$\overline{\Bbb F}_{q}(B)= \overline{\Bbb F}_{q}(t)$ with respect to $t$.
Since
\begin{eqnarray*}
  \det
     \left(
       \begin{array}{c}
         \varphi \\
         D_t^{(1)}\varphi \\
         D_t^{(\nu)}\varphi
       \end{array}
     \right)
     &=& 
     \det
     \left(
        \begin{array}{ccc}
         1 & t+t^q & t^{q+1} \\
         0 &   1   & t^q \\
         0 & \binom{q}{\nu}t^{q-\nu}&\binom{q+1}{\nu}t^{q+1-\nu}
        \end{array}
     \right)\\
     &=&
     \left\{
       \begin{array}{cl}
         0 & (q > \nu >1) \\
         t-t^q& (\nu =q)
       \end{array}
     \right.,
\end{eqnarray*}
the order-sequence for $\Phi$ is $\{0, 1, q\}$.

Since the tangent line at $\Phi(P_{\alpha})$ is given by (\ref{tangent}),
$P_{\alpha}^q = (1,(\alpha^q+ \alpha)^q, (\alpha^{q+1})^q)$
lies on it.
Hence $B$ is $q$-Frobenius nonclassical.
The set of inflection points of $B$ is given by the support of
the Wronskian divisor:
\[
(1+2)(q+1)P_{\infty} + 
{\rm div}\left( \det \left(\begin{array}{c}
         \varphi \\
         D_t^{(1)}\varphi \\
         D_t^{(q)}\varphi
       \end{array} \right) \right)
       +(0+1+q){\rm div}\, dt
     = \sum_{P \in {\Bbb P}^1({\Bbb F}_q)}\, P.
\]
Hence ${\rm Flex}\, B = \Phi({\Bbb P}^1({\Bbb F}_q))$,
which, together with (\ref{ratptB}),
implies (vi).
\qed

\gyokan

As was mentioned in Introduction, $N_q(B)$ attains
the Aubry-Perret bound for rational curves.
\begin{corollary}\label{NumberRatPtB}
The Ballico-Hefez curve $B$ is of degree $q+1$, and
$N_q(B) = q+1 + \frac{q(q-1)}{2}.$
\end{corollary}
\begin{corollary}\label{zetaB}
The zeta function of $B$ is
$
Z_B(T) = \frac{(1+T)^{\frac{q^2-q}{2}}}{(1-T)(1-qT)}.
$
\end{corollary}
\proof
Since $\Phi: {\Bbb P}^1 \to B$ is the normalization of $B$,
we can calculate the zeta function of $B$ by \cite[Th. 2.1]{aub-per}
together with the following informations:
${\rm Sing}\, B$ consists of $\frac{q^2-q}{2}$ points
that are ${\Bbb F}_q$-rational, and $\Phi^{-1}(\Phi(P)) = \{ P, P^q \}$
with $P \in {\Bbb P}^1({\Bbb F}_{q^2}) \setminus {\Bbb P}^1({\Bbb F}_q)$.
\qed

\begin{remark}\label{otherSim}
Corollaries \ref{NumberRatPtB} and \ref{zetaB} suggest another similarity
between the Ballico-Hefez curve and the Hermitian curve.
In general, if $C'$ is an irreducible curve over ${\Bbb F}_q$
with the normalization ${\Bbb P}^1 \stackrel{\pi}{\to} C'$,
then the zeta function of $C'$ is of the form
$\frac{L_{C'}(T)}{(1-T)(1-qT)}$
where $L_{C'}(T)$ is a polynomial of degree
${}^{\#}\pi^{-1}({\rm Sing}\, C') - {}^{\#}{\rm Sing}\, C'$
(say, $\Delta_{C'}$).
Let $\{ \beta_1, \ldots , \beta_{\Delta_{C'}} \}$ be the set of
reciprocal roots of $L_{C'}$.
Then $|\beta_i|=1$ ($i=1, \ldots , \Delta_{C'}$) and
$N_{q^r}(C')=q^r+1 - \sum_{i=1}^{\Delta_{C'}} \beta_i^{r}$
(see, \cite{aub-per}).
Especially, $N_q(C') \leq q+1+\Delta_{C'}$, and
if equality holds, then we have
$\Delta_{C'} \leq \frac{q^2-q}{2}$ by Ihara's argument \cite{iha}.
Actually equality holds in both inequalities for the Ballico-Hefez curve.
Comparing with R\"{u}ck and Stichtenoth's characterization
of Hermitian curves \cite{ruc-sti}, we might expect Ballico-Hefez curves
to be characterized among rational curves by those two properties.
\end{remark}

\gyokan

The following lemma will be used later.
\begin{lemma}\label{tangentSmooth}
For each ${\Bbb F}_q$-point $Q$ of ${\Bbb P}^2$ which does not lie on $B$,
there are two points $P_1, P_2 \in B({\Bbb F}_q) \setminus {\rm Sing}\, B$
such that $T_{P_1}B \cap T_{P_2}B = \{Q\}$.
Moreover the pair $\{ P_1, P_2\}$ is uniquely determined by $Q$.
\end{lemma}
\proof
Since $i(B.T_{P_i}B;P_i) = q+1$ by Theorem~\ref{theorem1},
$T_{P_i}B \cap B =\{ P_i\}$.
Hence the map
\[
  \begin{array}{ccc}
   S^2(B({\Bbb F}_q)\setminus {\rm Sing}\, B) \setminus \Delta &
    \to & {\Bbb P}^2({\Bbb F}_q) \setminus B({\Bbb F}_q) \\
    \{ P_1, P_2 \} & \mapsto & T_{P_1}B \cap T_{P_2}B
  \end{array}
\]
is well-defined, where $S^2(B({\Bbb F}_q)\setminus {\rm Sing}\, B)$
denotes the symmetric product of $B({\Bbb F}_q)\setminus {\rm Sing}\, B$
and $\Delta$ its diagonal subset.
Since the source and the target of this map have the same cardinality
$\frac{q(q+1)}{2}$,
it is enough to show the following fact;
{\em
Let $P_i' = (\alpha_i, \beta_i) \in {\Bbb P}^1$ $(i = 1,2,3)$
be three distinct points. Then three embedded tangent lines
$d\Phi(T_{P_i', {\Bbb P}^1})$ $(i = 1,2,3)$ are not concurrent.
}
In fact, since $d\Phi(T_{P_i', {\Bbb P}^1})$ is spanned by
two vectors $(\alpha_i^q, \beta_i^q, 0)$ and 
$(0, \alpha_i^q, \beta_i^q)$ (see (\ref{dPhi})),
its equation is
$\beta_i^{2q}X - \alpha_i^q\beta_i^q Y + \alpha_i^{2q}Z=0$.
Since
\[
 \left|
    \begin{array}{ccc}
        \beta_1^{2q} & - \alpha_1^q\beta_1^q & \alpha_1^{2q}\\
        \beta_2^{2q} & - \alpha_2^q\beta_2^q & \alpha_2^{2q}\\
        \beta_3^{2q} & - \alpha_3^q\beta_3^q & \alpha_3^{2q}
    \end{array}
 \right|
 = -\prod_{i<j} (\alpha_i\beta_j - \alpha_j\beta_i)^q,
\]
those three lines are not concurrent.
\qed
\section{Geometry of $B$ with ${\Bbb F}_q$-lines, for $q$ odd}
The projective plane of ${\Bbb F}_q$-lines in ${\Bbb P}^2$
is denoted by $\check{\Bbb P}^2({\Bbb F}_q)$.
An ${\Bbb F}_q$-line $l$ is an $i$-line
if ${}^{\#}\left(B({\Bbb F}_q)\cap l \right) = i$.
The cardinality of
$\{ l \in \check{\Bbb P}^2({\Bbb F}_q) | \mbox{\rm $l$ is an $i$-line} \}$
is denoted by $\tau_i$.
Since $\deg B = q+1$, only $q+2$ numbers $\{ \tau_i | 0 \leq i \leq q+1 \}$
make sense.
The purpose of this and next sections is to determine the exact values
of the $\tau_i$'s.

In this section, we assume $q$ is odd.
\begin{lemma}\label{withConic}
Let $C_B$ be the conic defined by $4XZ-Y^2=0$.
Then
$B \cap C_B =B({\Bbb F}_q) \cap C_B({\Bbb F}_q)
= B({\Bbb F}_q) \setminus {\rm Sing}\, B$,
and $T_P(B) =T_P(C_B)$ for any $P \in B\cap C_B$.
\end{lemma}
\proof
Let $\Phi (\alpha, \beta) = P \in B({\Bbb F}_q) \setminus {\rm Sing}\, B$.
Then $(\alpha, \beta) \in {\Bbb P}^1({\Bbb F}_q)$.
Hence we may suppose $\alpha, \beta \in  {\Bbb F}_q$, and have
$P=(\alpha^{q+1}, \alpha^q \beta + \alpha \beta^q, \beta^{q+1})
    = (\alpha^2, 2\alpha\beta, \beta^2)$,
which lies on $C_B$.
So
\[
B({\Bbb F}_q) \setminus {\rm Sing}\, B
\subseteq C_B({\Bbb F}_q) \cap B({\Bbb F}_q)
\subseteq B \cap C_B.
\]
Let ${\bf b}_1(P)=(\alpha^q, \beta^q, 0)$,
 ${\bf b}_2(P)= (0, \alpha^q, \beta^q)$,
 ${\bf c}_1(P) = (2\alpha, 2\beta, 0)$, and
 ${\bf c}_2(P) = (0, 2\alpha, 2\beta)$.
 Then $T_PB$ is spanned by ${\bf b}_1(P)$ and ${\bf b}_2(P)$, and
 $T_PC_B$ by ${\bf c}_1(P)$ and ${\bf c}_2(P)$.
Since $\alpha, \beta \in {\Bbb F}_q$,
$T_PB = T_PC_B$.
Hence
\[
2(q+1) = (B.C_B) =
\sum_{P \in B({\Bbb F}_q) \setminus {\rm Sing}\, B}
i(B.C_B; P) \geq 2(q+1),
\]
which means 
$B \cap C_B = B({\Bbb F}_q) \setminus {\rm Sing}\, B$.
\qed

Now we give geometric characterization of the set of rational points
$B({\Bbb F}_q)$ of the Ballico-Hefez curve for odd $q$.
The configuration of $B({\Bbb F}_q)$ is more or less known.
We define two subsets ${\cal S}$ and ${\cal T}$ in
${\Bbb P}^2({\Bbb F}_q)$.
\begin{definition}\label{SandT}
\begin{enumerate}[(1)]
  \item Let $l_1, \ldots , l_{q+1}$ be $q+1$ ${\Bbb F}_q$-lines
  that form an arc in $\check{\Bbb P}^2({\Bbb F}_q)$,
  that is, no three of the $q+1$ lines are concurrent.
  ${\cal S}$ is the set
  ${\Bbb P}^2({\Bbb F}_q) \setminus \bigcup_{i<j} l_i\cap l_j$.
  \item Let $D$ be a conic\footnote{In our context, the word ``conic" connotes
that it is absolutely irreducible.} over ${\Bbb F}_q$.
  ${\cal T}$ denotes the internal points of $D({\Bbb F}_q)$
  together with $D({\Bbb F}_q)$, which 
  appeared in \cite[Example 12.6 (3)]{hir}.
\end{enumerate}
\end{definition}
\begin{theorem}\label{three_configurations}
The three subsets $B({\Bbb F}_q)$, ${\cal S}$ and ${\cal T}$ in
${\Bbb P}^2({\Bbb F}_q)$ are projectively equivalent over
${\Bbb F}_q$.
\end{theorem}
\proof
First we show that ${\cal S}$ can be constructed by the same way
of constructing ${\cal T}$.
By a theorem of Segre \cite{seg1955},
$l_1, \ldots , l_{q+1}$ lie on a conic $D'$ in
$\check{\Bbb P}^2({\Bbb F}_q)$.
Then the dual $D$ of $D'$ in ${\Bbb P}^2({\Bbb F}_q)$ is also a conic,
and each $l_i$ tangents to $D$.
Hence the point set
$\{ l_i \cap l_j | 1 \leq i <j \leq q+1 \}$
is the external points of $D({\Bbb F}_q)$.

Next we consider the case for $B({\Bbb F}_q)$ and ${\cal T}$.
By Lemmas \ref{tangentSmooth} and \ref{withConic},
${\Bbb P}^2({\Bbb F}_q) \setminus B({\Bbb F}_q)$
is the external points of $C_B({\Bbb F}_q)$.
Therefore $B({\Bbb F}_q)$ coincides with ${\cal T}$
made from $C_B$.
Since any two conics over ${\Bbb F}_q$ are projectively equivalent,
so are those three sets.
\qed

\begin{corollary}\label{characterization_i_line_odd}
 \begin{enumerate}[{\rm (i)}]
   \item If $\tau_i \neq 0$, then $i=1$ or $\frac{q+1}{2}$ or $\frac{q+3}{2}$.
   \item $\tau_1 = q+1${\rm ;} $\tau_{\frac{q+1}{2}} = \frac{q(q-1)}{2}${\rm ;}
    $\tau_{\frac{q+3}{2}} = \frac{q(q+1)}{2}$.
 \end{enumerate}
\end{corollary}
\proof
We count these numbers by using the configuration ${\cal T}$ in
Definition~\ref{SandT}.
Let $l$ be an ${\Bbb F}_q$-line which tangents to the conic $D$
at $P \in D({\Bbb F}_q)$.
For each point $Q \in l({\Bbb F}_q) \setminus \{P\}$,
there exists another tangent ${\Bbb F}_q$-line passing through $Q$,
that is, $Q$ is an external point of $D({\Bbb F}_q)$.
Therefore any tangent line is $1$-line.
For an ${\Bbb F}_q$-line which does not tangent to $D$
at any ${\Bbb F}_q$-points, 
there are exactly $\frac{q+1}{2}$ external points
if the line does not meet $D({\Bbb F}_q)$, or
$\frac{q-1}{2}$ if it meets $D({\Bbb F}_q)$,
because there are exactly two ${\Bbb F}_q$-tangent lines passing through
an external point of $D({\Bbb F}_q)$.
Hence such a line is a $\frac{q+1}{2}$-line or a $\frac{q+3}{2}$-line.

For (ii), $\tau_1={}^{\#}D({\Bbb F}_q) =q+1$;
$\tau_{\frac{q+3}{2}}$ is the same as the number of ${\Bbb F}_q$-lines
joining two distinct points of $D({\Bbb F}_q)$,
which is $\frac{q(q+1)}{2}$.
Hence
$\tau_{\frac{q+1}{2}} = q^2+q+1 - \tau_1-\tau_{\frac{q+3}{2}}=\frac{q(q-1)}{2}$.\qed

\begin{remark}\label{characterization_i_line_T}
The proof of Corollary~\ref{characterization_i_line_odd} tells us a characterization of $i$-lines ($i=1, \frac{q+1}{2}, \frac{q+3}{2}$) for ${\cal T}$.
Let $l$ be an ${\Bbb F}_q$-line of ${\Bbb P}^2$.
Then $l$ is a $1$-line for ${\cal T}$ if and only if ${}^{\#}(l \cap D({\Bbb F}_q)) =1$;
a $\frac{q+1}{2}$-line for ${\cal T}$ if and only if ${}^{\#}(l \cap D({\Bbb F}_q)) =0$;
a $\frac{q+3}{2}$-line for ${\cal T}$ if and only if ${}^{\#}(l \cap D({\Bbb F}_q)) =2$.
\end{remark}
\section{Geometry of $B$ with ${\Bbb F}_q$-lines, for $q$ even}
The aim of this section is to determine the $\tau_i$'s for $B({\Bbb F}_q)$
with $q$ even,
which corresponds to Corollary~\ref{characterization_i_line_odd} for odd $q$.

\begin{proposition}\label{characterization_i_line_even}
Suppose $q$ is a power of $2$.
\begin{enumerate}[{\rm (i)}]
 \item The $q+1$ points of $B({\Bbb F}_q) \setminus {\rm Sing}\,  B$
are collinear.
 \item If $\tau_i \neq 0$, then $i=1$ or $\frac{q}{2}+1$ or $q+1$.
 \item $\tau_1 = q+1${\rm ;} $\tau_{\frac{q}{2}+1} = q^2-1${\rm ;}
    $\tau_{q+1} = 1$.
\end{enumerate}
\end{proposition}
\proof
Let $\Phi (\alpha, \beta) = P \in B({\Bbb F}_q) \setminus {\rm Sing}\, B$.
Since $(\alpha, \beta) \in {\Bbb P}^1({\Bbb F}_q)$,
$P=(\alpha^2, 0, \beta^2)$,
that is, $P$ lies on the line $Y=0$.
Let $l_0 = \{ Y=0 \}$.
Obviously $l_0({\Bbb F}_q) = B({\Bbb F}_q) \setminus {\rm Sing}\, B$
and $l_0$ is a $(q+1)$-line.
For any ${\Bbb F}_q$-line $l \neq l_0$,
put $n(l) = {}^{\#}( l \cap {\rm Sing}\, B)$.
Since any tangent line to a branch at a singular point $P$ of $B$ is not
${\Bbb F}_q$-line,
$i(B.l;P)=2$ for any ${\Bbb F}_q$-line $l$ passing through $P$.
Hence $1 + 2n(l) \leq q+1$.
Fix a point $P_0 \in l_0({\Bbb F}_q)$, and count the number of points of
${\rm Sing}\, B$ by using ${\Bbb F}_q$-lines passing through $P_0$.
Note that $l_0$ and $T_{P_0}B$ never meet ${\rm Sing}\, B$.
So we have
\[
\frac{q^2-q}{2} = 
\sum_{\scriptstyle l \ni P_0 \atop{\scriptstyle {\rm with} 
   \atop\scriptstyle l \neq l_0, T_{P_0}B}} n(l)
\leq \frac{q}{2}(q-1).
\]
Therefore an ${\Bbb F}_q$-line $l$ which neither $l_0$ nor the tangent line
to $B$ at an ${\Bbb F}_q$-point of $l_0$ is
a $(\frac{q}{2}+1)$-line,
and the number of such lines is $(q+1)(q-1)$.
Obviously the tangent line to $B$ at an ${\Bbb F}_q$-point of $l_0$ is
a $1$-line, and the number of such lines is $q+1$.
This completes the proof.
\qed

\section{Codes from Ballico-Hefez curves}
\begin{explanation}
{\em Codes from a subset of ${\Bbb P}^2({\Bbb F}_q)$.}
Let $S$ be a subset of ${\Bbb P}^2({\Bbb F}_q)$ which consists of $s$ elements.
For each point $P \in S$, we fix a representative $(a_0, a_1, a_2)$
of its coordinates with $a_0, a_1, a_2 \in {\Bbb F}_q$.
Then for any homogeneous polynomial $F(X_0, X_1, X_2)$ over ${\Bbb F}_q$,
the value $F(P) \in {\Bbb F}_q$ is determined without ambiguity.

Let $\Gamma ({\cal O}(i))$ be the vector space of homogeneous polynomials over
${\Bbb F}_q$ of degree $i$.
Then the image of ${\Bbb F}_q$-linear map
\[
\Gamma ({\cal O}(i)) \ni F \mapsto (F(P))_{P \in S} \in ({\Bbb F}_q)^{s}
\]
gives a linear codes, which is denoted by $C_L(S, {\cal O}(i))$.
\end{explanation}

\begin{proposition}
 Suppose that $q$ is odd.
\begin{enumerate}[{\rm (i)}]
 \item $C_L(B({\Bbb F}_q), {\cal O}(1))$ is
a $\left[ \frac{q^2+q+2}{2}, 3, \frac{q^2-1}{2} \right]_q$-code,
and it achieves the Griesmer bound.
 \item For $q \geq 5$, $C_L(B({\Bbb F}_q), {\cal O}(2))$ is
a $\left[ \frac{q^2+q+2}{2}, 6, \frac{q^2-q-4}{2} \right]_q$-code.
 \item  For $q \geq 7$, $C_L(B({\Bbb F}_q), {\cal O}(3))$ is
a $\left[ \frac{q^2+q+2}{2}, 10, \frac{q^2-2q-7}{2} \right]_q$-code.
\end{enumerate}
\end{proposition}
\proof
(i) The parameters are known by Corollaries~\ref{NumberRatPtB} and
\ref{characterization_i_line_odd}.
Since $\frac{q^2-1}{2} = \frac{q-1}{2} \cdot q + \frac{q-1}{2}$,
\[
\sum_{i=0}^{2} \left\lceil \left(\frac{q^2-1}{2}\right)/q^i \right\rceil
 = \frac{q^2-1}{2} + \left(\frac{q-1}{2} + 1\right) +1 = \frac{q^2+q+2}{2},
\]
which means the triple of parameters
$\left[ \frac{q^2+q+2}{2}, 3, \frac{q^2-1}{2} \right]_q$
achieves the Griesmer bound.

(ii) Let $D$ be a curve over ${\Bbb F}_q$ of degree $2$ in ${\Bbb P}^2$.
If $D$ is absolutely irreducible, then
${}^{\#}D({\Bbb F}_q) = q+1$;
if it is irreducible over ${\Bbb F}_q$ but not absolutely,
then ${}^{\#}D({\Bbb F}_q) = 1$.
Hence ${}^{\#}(B({\Bbb F}_q) \cap D) \leq  q+1$
for those two cases.
If $D$ is reducible over ${\Bbb F}_q$,
then ${}^{\#}(B({\Bbb F}_q) \cap D) \leq  q+3$ by Corollary~\ref{characterization_i_line_odd}.
In particular, no degree-two-curve contains $B({\Bbb F}_q)$,
and $C_L(B({\Bbb F}_q), {\cal O}(2))$ is
a $\left[ \frac{q^2+q+2}{2}, 6,  \geq \frac{q^2-q-4}{2} \right]_q$-code,
where $\geq \frac{q^2-q-4}{2}$ at the third parameter
means the minimum distance of this code is at least $\frac{q^2-q-4}{2}$.
For an external point $Q$ of $C_B({\Bbb F}_q)$,
which does not lie on $B({\Bbb F}_q)$ by Theorem~\ref{three_configurations},
there are exactly two ${\Bbb F}_q$-lines passing through $Q$
each of which tangents to $C_B$.
Hence there are exactly $\frac{q-1}{2}$ ${\Bbb F}_q$-lines passing through $Q$
each of which meets $C_B({\Bbb F}_q)$ at two points.
Since $q \geq 5$, we can choose two lines from those $\frac{q-1}{2}$ lines,
which are $\frac{q+3}{2}$-lines by Remark~\ref{characterization_i_line_T}.
The union of those two lines gives a codeword of weight $\frac{q^2-q-4}{2}$.

(iii) For $C_L(B({\Bbb F}_q), {\cal O}(3))$,
in order to verify that its parameters are
$\left[ \frac{q^2+q+2}{2}, 10, \geq \frac{q^2-2q-7}{2} \right]_q$,
it is enough to see that
\begin{equation}\label{claim_degree3}
{}^{\#}(B({\Bbb F}_q) \cap D) \leq \frac{3}{2}( q+3)
\end{equation}
for any curve $D$ over ${\Bbb F}_q$ of degree $3$.
If $D$ is absolutely irreducible,
then
${}^{\#}(B({\Bbb F}_q) \cap D) \leq q+1 + 2 \sqrt{q} \leq \frac{3}{2}(q+3)$.
In other cases, one can verify (\ref{claim_degree3}) easily
by using Corollary~\ref{characterization_i_line_odd}.
Equality in (\ref{claim_degree3}) holds if one takes three internal lines
of $C_B({\Bbb F}_q)$ passing through an assigned external point,
which is possible because $q \geq 7$.
\qed

When $q$ is a power of $2$,
$B({\Bbb F}_q)$ contains $q+1$ collinear points,
so we can't expect to obtain good codes from the total set.
\begin{proposition}
Let $q=2^e$ with $e>2$.
Then $C_L({\rm Sing}\, B, {\cal O}(1))$ is
a $\left[ \frac{q^2-q}{2}, 3, \frac{q^2-2q}{2} \right]_q$-code,
which achieves the Griesmer bound.
\end{proposition}
\proof
From Proposition~\ref{characterization_i_line_even} and its proof,
an $i$-line for ${\rm Sing}\, B$ exists if and only if $i=0$ or $\frac{q}{2}$.
Hence we have the first half of the assertion.
Since $\frac{q^2-2q}{2} = \frac{q-2}{2}\cdot q$,
we have the additional assertion.
\qed

\section{Generalization of the curve $B$}
We generalize the Ballico-Hefez curve in ${\Bbb P}^2$ to
rational curves in higher dimensional projective spaces. 
We will discuss only on the number of rational points on them.

\begin{notation}
Let $X_0, \dots , X_{n-1}$ be variables.
We consider $n+1$ elementary symmetric polynomials:
\[
\sigma_k (X_0, \dots , X_{n-1}):=
\sum_{i_1 < \dots < i_k}  X_{i_1}\cdots X_{i_k}
\text{ \ ($k= 0, \dots , n$),}
\]
where we understand $\sigma_0 =1$.
We also consider ``homogeneous" elementary symmetric polynomials
of degree $n$:
\begin{eqnarray*}
\lefteqn{
\Tilde{\sigma}_k(X_0, \dots , X_{n-1}; Y_0, \dots , Y_{n-1}):=}\\
& & \sum_{\scriptstyle i_1 < \dots < i_k; j_1 < \dots < j_{n-k} 
     \atop{\scriptstyle
      \text{with} \atop
     \scriptstyle
      \{i_1, \dots, i_k, j_1, \dots , j_{n-k} \} 
      =\{0, \dots , n-1\}} }
           X_{i_1}\cdots X_{i_k}Y_{j_1}\cdots Y_{n-k}.
\end{eqnarray*}
\end{notation}

\begin{definition}
We define $n+1$ homogeneous polynomials of degree
$q^{n-1}+q^{n-2}+\dots +1$ in $s$ and $t$ by
\[
\Tilde{\varphi}_k(s,t) :=
\Tilde{\sigma}_k(s, s^q, \dots, s^{q^{n-1}};t, t^q, \dots, t^{q^{n-1}}),
\]
and inhomogeneous ones by
\[
\varphi_k(s) := \Tilde{\varphi}_k(s,1) = \sigma_k(s, s^q, \dots, s^{q^{n-1}})
= \sum_{i_1 < \dots < i_k}  s^{q^{i_1}}\cdots s^{q^{i_k}}
\]
for $k = 0, 1, \dots , n$.

$B_n$ denotes a curve in ${\Bbb P}^n$ over ${\Bbb F}_q$
defined by the image of
\[
\Phi_n: {\Bbb P}^1 \ni (s,t)
 \mapsto (\Tilde{\varphi}_0(s,t), \Tilde{\varphi}_1(s,t),
 \dots , \Tilde{\varphi}_n(s,t) ) \in {\Bbb P}^n.
\]
\end{definition}

Since $B_2=B$, the curve $B_n$ is a generalization of
the Ballico-Hefez curve.

\begin{remark}
$B_n$ is nondegenerate in ${\Bbb P}^n$,
that is, there is no hyperplane of ${\Bbb P}^n$
containing $B_n$.

In fact, suppose
$\sum_{i=0}^n \alpha_i \varphi_i(s) =0$
for $\alpha_0, \dots, \alpha_n \in \overline{\Bbb F}_q$.
Let $D_s^{(\nu)}$ be the $\nu$-th Hasse derivation with respect to $s$.
Note that for $k$ and $i_0 < \dots <i_l$,
$D_s^{(1+q+q^2+ \dots +q^k)}s^{q^{i_0}+q^{i_1}+\dots+q^{i_l}} \neq 0$
if and only if $l \geq k$ and $i_{\mu} =\mu$ ($\mu = 0, \dots , l$),
and in this case,
$
D_s^{(1+q+q^2+ \dots +q^k)}s^{q^{i_0}+q^{i_1}+\dots+q^{i_l}}
= s^{q^{i_{k+1}}+\dots+q^{i_l}}.
$
Therefore
\[
0 = D_s^{(1+q+q^2+ \dots +q^k)}
   \left(\sum_{i=0}^n \alpha_i \varphi_i(s) \right)
   |_{s=0} =\alpha_k,
\]
which means that no hyperplane contains $B_n$.
\end{remark}
\begin{theorem}\label{generalization}
\begin{enumerate}[{\rm (i)}]
\item For $P \in {\Bbb P}^1$, the map arising from $\Phi_n$
on tangent spaces
\[
d\Phi_{n, P}:  T_{P, {\Bbb P}^1} \to T_{\Phi_n(P), {\Bbb P}^n}
\]
is injective.
\item For $P \in {\Bbb P}^1$, ${}^{\#} \Phi_n^{-1}(\Phi_n(P))> 1 $
if and only if
$P \in {\Bbb P}^1({\Bbb F}_{q^n}) \setminus {\Bbb P}^1({\Bbb F}_{q})$.
In this case, if ${\Bbb F}_{q}(P) = {\Bbb F}_{q^k}$,
then
\[
\Phi_n^{-1}(\Phi_n(P))
= \{ P, P^q, \dots P^{q^{k-1}}\},
\]
and
\[
\Phi_n(P)=\Phi_n(P^q)= \dots =\Phi_n(P^{q^{k-1}}) 
    \in {\Bbb P}^n({\Bbb F}_{q}).
\]
\item $B_n({\Bbb F}_{q}) = \Phi_n({\Bbb P}^1({\Bbb F}_{q^n}))$.
\end{enumerate}
\end{theorem}
\proof
(i) Since
\[
\left(
 \begin{array}{c}
  \frac{\partial \Phi_n}{\partial s} \\
  \frac{\partial \Phi_n}{\partial t}
 \end{array}
\right)
=
\left(
  \begin{array}{ccccccc}
   0   & t^{q+ \dots +q^{n-1}}& \ast & \cdots &\ast 
                           &\ast &s^{q+ \dots +q^{n-1}}\\
t^{q+ \dots +q^{n-1}}&\ast    & \ast & \cdots &\ast
                           &s^{q+ \dots +q^{n-1}}&0
  \end{array}
\right),
\]
the rank of 
$
\left(
 \begin{array}{c}
  \frac{\partial \Phi_n}{\partial s} \\
  \frac{\partial \Phi_n}{\partial t}
 \end{array}
\right)
$
is $2$ for any point $(s,t) \in {\Bbb P}^1$.

(ii) We may suppose $t=1$.
For two points $P_{\alpha}=(\alpha, 1), P_{\beta}=(\beta,1) \in {\Bbb P}^1$,
$\Phi_n(P_{\alpha}) = \Phi_n(P_{\beta})$
if and only if $\varphi_k(\alpha) = \varphi_k(\beta)$
for $k=0,1,\dots,n$.
Since
$
\prod_{i=0}^{n-1} (X - s^{q^i}) 
     = \sum_{k=0}^{n} (-1)^k \varphi_k(s)X^{n-k},
$
those conditions are equivalent to the condition
\[
\{\alpha, \alpha^q, \dots, \alpha^{q^{n-1}} \}=
    \{\beta, \beta^q, \dots, \beta^{q^{n-1}}\}
\]
with counting multiplicity.
Here we need the following lemma, which completes the proof of (ii).
\begin{lemma}
Let two elements $\alpha$ and $\beta$ of $\overline{\Bbb F}_q$
be distinct from each other.
If
$
\{\alpha, \alpha^q, \dots, \alpha^{q^{n-1}} \}=
    \{\beta, \beta^q, \dots, \beta^{q^{n-1}}\}
$
with counting multiplicity,
then ${\Bbb F}_q(\alpha) = {\Bbb F}_q(\beta)$
is a subfield of ${\Bbb F}_{q^n}$.
\end{lemma}
\proof
Let ${\Bbb F}_q(\alpha)= {\Bbb F}_{q^k}$.
Obviously ${\Bbb F}_q(\beta)= {\Bbb F}_{q^k}$ also.
Since $\beta = \alpha^{q^i}$ for some $i$ with $1 \leq i \leq n-1$,
$\alpha^{q^n}=(\alpha^{q^i})^{q^{n-i}} = \beta^{q^{n-i}}
 \in \{ \alpha, \alpha^q, \dots , \alpha^{q^{n-1}}\}$.
Hence $\alpha^{q^n} = \alpha^{q^j}$ for some $j$ with $0 \leq j \leq n-1$.
Hence $k|n-j$, particularly $k \leq n$.
Let $n = uk +r$ with $0 \leq r <k$.
Suppose that $r >0$.
Then each of the $r$ elements
$\{ \alpha, \alpha^{q}, \dots, \alpha^{q^{r-1}} \}$
appears $u+1$ times in the total set
$\{\alpha, \alpha^q, \dots, \alpha^{q^{n-1}} \}$,
and each of the $k-r$ elements
$\{ \alpha^{q^r}, \dots, \alpha^{q^{k-1}} \}$
$u$-times.
On the other hand, since $\beta = \alpha^{q^i}$,
each of the $r$ elements
$\{ \alpha^{q^i}, \alpha^{q^{i+1}}, \dots, \alpha^{q^{i+r-1}} \}$
appears $u+1$ times in the total set, and
each of the $k-r$ elements
$\{ \alpha^{q^{i+r}}, \dots, \alpha^{q^{i+k-1}} \}$
$u$-times.
Therefore both $\alpha^{q^{r-1}}$ and $\alpha^{q^{i+r-1}}$
appear $u+1$ times in the total set, but their $q$-th power $u$ times.
Hence those elements must coincide.
Hence $\alpha = \alpha^{q^i}=\beta$, which is a contradiction.
\qed

\gyokan

\noindent
{\em Continuation of the proof of Theorem~{\rm \ref{generalization}.}}
(iii) It is obvious that
$\Phi_n(P_{\alpha}) \in B({\Bbb F}_q)$
for any $\alpha \in {\Bbb F}_{q^n}$.
Conversely, if $\Phi_n(P_{\alpha}) \in B({\Bbb F}_q)$,
then
$\varphi_n(\alpha) = \alpha^{1+q+ \dots +q^{n-1}} \in {\Bbb F}_q$.
Hence
$1=(\alpha^{1+q+ \dots +q^{n-1}})^{q-1}=\alpha^{q^n-1}$ if $\alpha \neq 0$,
which means $\alpha \in {\Bbb F}_{q^n}$.
\qed

\gyokan

From Theorem~\ref{generalization},
we know the number $N_q(B_n)$ of ${\Bbb F}_q$-points of $B_n$.
\begin{theorem}
\[
N_q(B_n)
= \frac{1}{n}q^n +
     \sum_{\scriptstyle d|n \atop{\scriptstyle\text{with}
               \atop\scriptstyle d\neq n }}
                 \frac{1}{d} \prod_{\scriptstyle l : \text{prime }
                                   \atop {\scriptstyle\text{with}
                                   \atop\scriptstyle
                                    l | \frac{n}{d}}}
                                 \left(1-\frac{1}{l}\right)q^d
                                                    +1.
\]
\end{theorem}
\proof
By (ii) and (iii) of Theorem~\ref{generalization},
\[
N_q(B_n) = \sum_{k|n} 
 \frac{1}{k} {}^{\#} \{P \in {\Bbb P}^1 \mid {\Bbb F}_q(P) ={\Bbb F}_{q^k} \}.
\]
Since ${\Bbb P}^1(\overline{\Bbb F}_q)=\overline{\Bbb F}_q \cup \{(1,0)\}$,
\[
{}^{\#} \{P \in {\Bbb P}^1 \mid {\Bbb F}_q(P) ={\Bbb F}_{q^k} \}
 = \left\{
    \begin{array}{lc}
    {}^{\#}\{\alpha  \in {\Bbb F}_{q^k} \mid 
        {\Bbb F}_q(\alpha) ={\Bbb F}_{q^k}  \}& \text{ if } k>1 \\
     {}^{\#}{\Bbb F}_q + 1 & \text{ if } k=1.
    \end{array}
   \right.
\]
Now we compute 
${}^{\#}\{\alpha  \in {\Bbb F}_{q^k} \mid 
        {\Bbb F}_q(\alpha) ={\Bbb F}_{q^k}  \}$ exactly.
Let $k = l_{k,1}^{e_1}\dots l_{k,r_k}^{e_{r_k}}$
be the decomposition of $k$ as a product of powers of distinct prime numbers.
Here we understand the right hand of the decomposition of $1$ to be empty.
Since
\[
\{\alpha  \in {\Bbb F}_{q^k} \mid 
        {\Bbb F}_q(\alpha) ={\Bbb F}_{q^k} \}
        = {\Bbb F}_{q^k}
            \setminus \bigcup_{i=1}^{r_k} {\Bbb F}_{q^{k/l_{k,i}}},
\]
the cardinality of this set is
\[
 q^k - \sum_{i=1}^{r_k} q^{k/l_{k,i}} + \sum_{i<j} q^{k/l_{k,i}l_{k,j}}
 - \cdots 
= \sum_{s=0}^{r_k} (-1)^s \sum_{i_1< \dots <i_s}
        q^{k/l_{k,i_1} \dots l_{k,i_s}}.
\]
When $k=1$, we understand the above sum to be just $q$.
Hence
\[
N_q(B_n) = \sum_{k|n}\frac{1}{k}\sum_{s=0}^{r_k} (-1)^s \sum_{i_1< \dots <i_s}
        q^{k/l_{k,i_1} \dots l_{k,i_s}} +1.
\]
For each $d$ with $d|n$, we gather together the terms in $q^d$.
If $d=n$, the only term is $\frac{1}{n}q^n$.
If $d<n$, put $\frac{n}{d} = m_1^{f_1}\cdots m_u^{f_u}$,
where $m_1, \dots , m_u$ are distinct prime numbers
with $f_1, \dots , f_u \geq 1$.
It is obvious that for a fixed $d$,
\[
\{ k \mid k/l_{k,j_1}\cdots l_{k, j_s} =d \text{ for some $j_1< \dots < j_s$} \}
= \{ m_{i_1}\cdots m_{i_s} d \mid i_1 < \dots <i_s \}.
\]
Hence the coefficient of $q^d$ is
\[
\sum_{s=0}^{u} \sum_{i_1 < \dots <i_s}(-1)^s
  \frac{1}{m_{i_1}\cdots m_{i_s}d}
  = \frac{1}{d}\prod_{i=1}^{u}\left(
                              1-\frac{1}{m_i}
                             \right)
 = \frac{1}{d}\prod_{\scriptstyle l|\frac{n}{d} \atop
                       \scriptstyle l: \text{ prime}}
                       \left( 1 - \frac{1}{l}\right).
\]
This completes the proof.

\gyokan

\begin{remark}
The order-sequence, the $q$-Frobenius index and the $q$-Frobenius
order-sequence for $\Phi_n: {\Bbb P}^1 \to B_n$ can be also calculated.
Namely, the order-sequence for $\Phi_n$ is
$
0<1< q< q^2 < \dots < q^{n-1},
$
the $q$-Frobenius index is $1$, and hence 
the $q$-Frobenius order-sequence is
$
0<q< q^2 < \dots < q^{n-1}.
$
At present, our proof of this fact is not concise.
So we will discuss it in a separate paper.
\end{remark}

\end{document}